\documentclass[10pt, reqno]{amsart}

\usepackage{amsmath,graphics}
\usepackage{amsfonts,amssymb}

\theoremstyle{plain}
\newtheorem*{theorem*}{Theorem}
\newtheorem*{lemma*} {Lemma}
\newtheorem*{corollary*} {Corollary}
\newtheorem*{proposition*} {Proposition}
\newtheorem{theorem}{Theorem}[section]
\newtheorem{lemma}[theorem]{Lemma}
\newtheorem{corollary}[theorem]{Corollary}
\newtheorem{proposition}[theorem]{Proposition}

\theoremstyle{definition}
\newtheorem{defn}[theorem]{Definition}

\def \R {\mathbb{R}}

\def \bn{\begin{enumerate}}
\def \en{\end{enumerate}}
\def \bdm{\begin{displaymath}}
\def \edm{\end{displaymath}}

\def \bp{\begin{proof}}
\def\ep{\end{proof}}
\def\be{\begin{equation}}
\def\ee{\end{equation}}

\begin{document}

\title{Cohomology of flat bundles and a Chern-Simons functional}
\author{Yanghyun Byun}
\address{Department of Mathematics, Hanyang University, Seoul, Korea}
\email{yhbyun@hanyang.ac.kr}
\author{Joohee Kim}
\address{Department of Mathematics, Hanyang University, Seoul, Korea}
\email{kjh0423@hanyang.ac.kr}
\subjclass[2010]{53C05, 55R20}

\begin{abstract}
We show that a flat principal bundle with compact connected structure group
and its adjoint bundles of Lie groups have the same cohomology
as the trivial bundle, which is done by proving
they satisfy the condition for the Leray-Hirsch theorem.
This information
has been used to construct a cohomology class of the adjoint bundle of a flat bundle
whose restriction to each fiber is the class of the Maurer-Cartan 3-form.
Then we use this result to define an invariant of a gauge
transformation of a flat bundle which describes the effect of the gauge transformation on a Chern-Simons
functional.
\end{abstract}
\maketitle

\tableofcontents

\section{Introduction}
The main themes of this paper are to calculate the cohomology of flat bundles,
to use the result for a rigourous definition of the degree of
a gauge transformation of a flat bundle and subsequently
to describe the effect of a gauge transformation on a Chern-Simons functional.
However a large portion
is devoted to introducing the relevant basic materials.
We hope this makes the manuscript truly self-contained in the sense that
a reader does not have to gather the references to learn or get reminded
of some basic notions in gauge theory. The contents of this paper are
the results of a few years of
collaboration of the authors and constitute the thesis submitted to Hanyang University
as a partial fulfillment for the Ph. D. of the second author.

Let $G$ be a Lie group with the Lie algebra $\mathcal{G}$.
We consider a principal $G$-bundle $\pi:P \rightarrow M$
over a smooth manifold $M$. Let $\mathcal{A}(P)\subset \Omega^1(P; \mathcal{G})$
denote the set of all connections on $P$. Choose an $A_0 \in \mathcal{A}(P)$ and  an
Ad-invariant inner product $\langle \cdot , \cdot \rangle$ on $\mathcal{G}$.
Let $A$ be any connection and
$F_A\in \Omega^2(P;\mathcal{G})$ denote
the curvature of $A$.
Then we consider the $3$-form
$$cs(A)=\langle (F_A+ F_{A_0})\wedge (A-A_0)\rangle - {\frac{1}{6}} \langle (A-A_0)\wedge [(A-A_0)\wedge(A-A_0)]\rangle.$$
The notational conventions used in this formula are explained in \S 3 below.

Since $cs(A)$ is horizontal and invariant, it may be regarded as a form on $M$.
Now assume that $M$ is an oriented closed $3$-manifold. Then we may define a function,
$CS:\mathcal{A}(P)\rightarrow \R$,
by the equation:
\bdm CS(A)=\int_M cs(A),\tag{1-1}\edm
which we refer to as the Chern-Simons functional of $A$ defined
by choosing the reference connection $A_0$.

The equality below reveals the Chern-Weil theory origin of the form $cs(A)$(see \S 5):
\bdm \langle F_A\wedge F_A \rangle - \langle F_{A_0}\wedge F_{A_0} \rangle
        = d (cs(A)).\tag{1-2} \edm

Note that $\mathcal{A}(P)$ is an Affine subspace of $\Omega^1(P;\mathcal{G})$
modeled after $\bar \Omega^1(P; \mathcal{G})$, the subspace of all horizontal equivariant
1-forms.  Then the following identity makes a good sense(see \S 5):
\bdm d CS_A(a) = 2 \int_M \langle F_A \wedge a\rangle,\tag{1-3} \edm
for any $A\in \mathcal{A}(P)$ and any $a \in T_A \mathcal{A}(P)\equiv \bar \Omega^1(P; \mathcal{G})$.
Thus the critical points of the Chern-Simons functional
are the flat connections.  Furthermore the RHS of (1-3) does not depend on the reference connection and therefore
we know that $CS_1 -CS_0$ is the constant $CS_1(A_0)=-CS_0(A_1)$ if $CS_0$ and $CS_1$ are
the Chern-Simons functionals given by reference connections $A_0$ and $A_1$ respectively.

The Chern-Simons functional given by (1-1) has the advantage that it can be
applied to any principal $G$-bundle $P$ over an oriented closed 3-manifold,
even when $P$ is not
trivializable. This form of functional has been used
for instance in \cite{Dostoglou}, \cite{Mrowka} and \cite{Wehrheim}.

Now let $\varphi: P \to P$ be a gauge transformation.
Then $CS(\varphi^*A)-CS(A)$ is
an invariant of the path component of $\varphi$ in the group
${\rm GA}(P)$ of gauge transformations, which depends
neither on $A$ nor on the reference connection, as observed
in \cite{Dostoglou}(see 5.2 below). Therefore $CS(\varphi^*A)-CS(A)$ is a well-defined
homotopy invariant of $\varphi$ and one may ask
whether this invariant can be described in some other way.
Answering this is one of the main results of the note.

The answer depends on the assumptions that the Lie group
$G$ is connected and compact and that $P$ is {\it flat}, that is, $P$ admits
a flat connection.
Under these assumptions we show that the equality
\bdm CS(\varphi^*A)-CS(A)= \deg \varphi \tag{1-4}\edm
holds for an invariant $\deg \varphi$ which is defined
independently of the LHS. Also we will prove that
$\deg \varphi$ is integer valued if an appropriate Ad-invariant inner product
$\langle \cdot , \cdot \rangle$ on $\mathcal{G}$ is chosen.

We use the assumptions as follows.
We first observe that a flat connection $A$ can be used to
define a homomorphism between the cohomology rings
$$E_A : H^*_{\rm dR}(G) \rightarrow H^*_{\rm dR}(Ad\  P),$$
if $G$ is connected and compact(see 6.5 below).
Here $Ad\ P=P\times_G G$ denotes the adjoint bundle of Lie groups.
Let $p\in P$ and $\iota_p: G \rightarrow Ad\ P$ be the map defined by
$\iota_p(g)=[p,g]$. Then in fact the homomorphism $E_A$ is
such that $\iota^*_pE_A$ is
the identity homomorphism on
$H^*_{\rm dR}(G)$. This means that the Leray-Hirsch theorem(cf.\ \cite{Hatcher})
holds in a stronger sense in that $E_A$ is not only a surjective homomorphism between vector spaces
such that $\iota^*_pE_A$ is
the identity but also a homomorphism between rings.
It follows that we have
$H^*_{\rm dR}(Ad\ P)$ is isomorphic to $H^*_{\rm dR}(M)\otimes H^*_{\rm dR}(G)$
not only as an $H^*_{\rm dR}(M)$-module but also as an $\mathbb{R}$-algebra.

Let $\theta\in \Omega^1(G;\mathcal{G})$ denote the Maurer-Cartan form(see \S 2).
Then we have a closed 3-form
$$\Theta =-\frac{1}{6}\langle \theta\wedge [\theta\wedge\theta]\rangle\in \Omega^3(G;\mathbb{R}),$$
which we call the Maurer-Cartan 3-form(see \S 3).
We will prove that $E_A[\Theta]\in H^3_{\rm dR}(Ad\ P)$ does not depend on the choice of the flat connection $A$(see 7.2 below).
We will also prove that $E_A[\Theta]$ is an integral class for
an appropriate choice of the Ad-invariant inner product $\langle \cdot , \cdot \rangle$.

So far no restriction on $M$ has been needed. Now assume $M$ is an oriented closed $3$-manifold.
Let $u: P\rightarrow G$ denote the map defined by $\varphi(p)=pu(p)$ for any $p\in P$ and
let $\hat u :M \rightarrow Ad P$ be defined by $\hat u(x) = [p, u(p)]$ for any $x\in M$
by choosing any $p\in P_x=\pi^{-1}\{ x\}$.
Then we define $\deg \varphi$ by
\bdm \deg \varphi=\int_M \hat u^*E_A[\Theta], \tag{1-5}\edm
choosing any flat connection $A\in \mathcal{A}(P)$(see 8.3 below).
Then (1-4) can be established
by a straightforward argument(see 8.4 below).

A definition of the degree
of a gauge transformation of a nontrivial bundle
has been given in case $G={\rm SO}(3)$ in \cite{Dostoglou},
which utilizes the language of algebraic topology
and therefore is quite different from
our definition given by (1-5)
which uses the language of differential category.
In fact a closed $3$-form
whose restriction to each fiber is the Maurer-Cartan 3-form appeared in a work(\cite{Wehrheim})
by K.\ Wehrheim. She exploited the 3-form
to define the degree of a gauge transformation, which we also followed in this note.
However she have justified neither the existence nor the uniqueness of such a 3-form.
From the view point provided by this paper
her argument as a whole is unharmed by this negligence
since all the bundles in concern in her work were flat.
In fact this paper appears the first to pinpoint flatness
as a condition under which
a principal bundle may admit such a 3-form.
The authors do not know
whether it exists for a general principal bundle.
They also do not know whether the condition that the base manifold $M$
is of dimension 3 guarantees such a 3-form should exist
in $H^3(Ad\ P; \mathbb{R})$ for a general bundle $P$ over $M$.

We also would like to mention the fact that a method to
define a Chern-Simons functional
on a general bundle over an oriented closed $3$-manifold with compact structure group
has been proposed by R.\ Dijkgraaf and E.\ Witten(\cite{Dijkgraaf})
(see also Appendix of \cite{Freed} by D.\ S.\ Freed).
The resulting functional
is $\mathbb{R}/\mathbb{Z}$-valued from the beginning.
Therefore in this method one simply observes that the functional
is invariant under gauge transformations and a need to define their degrees
does not surface explicitly. It appears that
the condition(s) under which the degree of a gauge transformation can be defined
has not yet been fully understood and the current note provides
one of them.

\section{Connections}
It is not easy to pinpoint a single appropriate
reference for this section. We believe that our definition of a connection
in terms of the Maurer-Cartan 1-form can be easily seen to be equivalent
to the more traditional one given in such as in \cite{bleeker}.
Our definition of the twisted derivative $d_A \alpha$ of a vector valued form $\alpha$
on $P$ by a connection $A$ is in fact that of \cite{bleeker} or of \cite{bleeker-booss}, adopting a slightly
different notation. We will also briefly review the Chern-Weil theory.

Let  $G$ be a Lie group. A surjective smooth map $\pi:P\rightarrow M $
 together with a smooth right action $P \times\ G\rightarrow P $ is a principal $G$-bundle
 if the following conditions are satisfied:

   1) $G$ acts freely and transitively on each fiber $ P_x = \pi^{-1}\ (x)$.

   2) Local triviality condition holds. That is, for each  $x \in\ M$,
   there is an open neighborhood $U$ of $x$ and a diffeomorphism $ h :\pi ^ {-1}\ (U)  \rightarrow  U \times\ G  $
   such that, writing $ h (p) =(\pi(p),\varphi(p))$ for any
 $p \in \pi ^ {-1}\ (U)$, we have $\varphi (pg)=\varphi(p)g$, for any $ g\in G$.

    Before we explain connections, we first introduce the Maurer-Cartan 1-form $\theta\in \Omega^1 (G;\mathcal{G})$,
    where $\mathcal{G}=T_eG$ is the Lie algebra of $G$.
    In fact $\theta$ is defined by $\theta (X)= dL_{g^{-1}}(X)\in\mathcal{G}$ for any $X\in T_g G ,\ g\in G $.
    Then we have:
    $$L_{g} ^* \theta = \theta ,\  R_g^*\theta =Ad_{g^{-1}}\circ \theta.$$
    The Maurer-Cartan 1-form will be needed also for the next section.

    Now choose any $p\in P_x$. Then let $\kappa_p : P_x \rightarrow G$ be the map defined by $\kappa_p (pg)=g$
    for any $ g\in G$.
    On the other hand let $\iota_x : P_x \rightarrow P$ be the inclusion.
    A connection $A$ is a $\mathcal{G}$-valued 1-form on $P$ which is equivariant, that is,
    $$R_g^* A =Ad_{g^{-1}}\circ A$$
    for any $g \in G $ and satisfies
    $$\kappa_p^* \theta =\iota_p ^* A \in \Omega^1 (P_x ; \mathcal{G}). $$

   The curvature $F_A \in\Omega ^2 (P;\mathcal{G})$ of $A$ is defined as the twisted derivative:
   $$F_A = d_A A.$$
   Then $F_A$ is equivarient, that is,
   $$R_g^*F_A =Ad_{g^{-1}}\circ F_A$$
   for any $g\in G$, and horizontal, that is,
   $$F_A(X,Y)=0$$
   if any of $X, Y\in T_pP$ is vertical.
   Also we have the Bianchi identity
   $$d_A F_A =0.$$

   Here we explain the twisted exterior derivative
   $d_A \alpha$ for any $\alpha \in \Omega ^k(P;V)$ where $V$ is any vector space.
   The horizontal subbundle $H$ of $TP$ determined by $A$ is the bundle whose fiber
   $H_p$ at $p$ is given by $H_p = Ker (A_p : T_pP \rightarrow \mathcal{G})$ for any $p\in P$.
   And the vertical subbundle $V$ of $TP$ is defined by $V_p = T_p P_{\pi(p)}$ for any $p\in P$.
   Then we have the decomposition $TP = H \oplus V$.
   Let $\pi_H : TP \rightarrow H$ be the associated projection.
   Then the twisted exterior derivative is defined by $d_A\alpha=(d\alpha)^H $ where
   $\beta ^H$, for any $\beta \in \Omega^k (P;V)$,
   is defined by
     $$\beta ^H (X_1,...,X_k)=\beta (\pi_H X_1,...,\pi_H X_k)$$
   for any $X_1,...,X_k \in T_pP,\ p\in P$.

   Now we review the Chern-Weil theory briefly(cf. \cite{Nomizu}).
   Let $I^k(\mathcal{G})$ be the set of all symmetric $k$-linear invariant functions,
     $f: \mathcal{G}\times ...\times\mathcal{G} \rightarrow \mathbb{R}$.
   That $f$ is invariant means that $f(Ad_g X_1 ,...,Ad_g X_k)= f(X_1,...,X_k)$
   for any $ g\in G $ and any $X_1,...,X_k \in \mathcal{G}$.
   For $f\in I^k(\mathcal{G})$ and $g\in I^l(\mathcal{G})$ we define
   $$(fg)(X_1,...,X_{(k+l)}) = \frac {1}{k!l!} \sum_{\sigma} f(X_{\sigma (1)},...,X_{\sigma (k)})g(X_{\sigma(k+1)},...,X_{\sigma(k+l)})$$
   where the summation is taken over all permutation $\sigma$ of $1,...,k+l$.
   Extending this multiplication to $I(\mathcal{G})=\bigoplus_{k=1}^{\infty}I^k(\mathcal{G})$,
   we make  $I(\mathcal{G})$ into a commutative associative algebra over $\mathbb{R}$.

   For each  $f\in I^k(\mathcal{G})$,  $f(F_A)$ is the real-valued $2k$-form on $P$ defined by
   $$f(F_A)(X_1,...,X_{2k}) = \frac{1}{2^k}\sum_{\sigma}\epsilon_{\sigma}f(F(X_{\sigma(1)},X_{\sigma(2)}),...,F(X_{\sigma(2k-1)},X_{\sigma(2k)}))$$
   for any $X_1,\cdots, X_{2k}\in T_pP$
   where $\epsilon_{\sigma}$ denotes the sign of the permutation $\sigma$.
   We note that $f(F_A)$ is horizontal and invariant therefore can be regarded as form on $M$. That is,
   there is a unique $\alpha\in \Omega^{2k}(M;\mathbb{R})$ such that $\pi^*\alpha=f(F_A)$ and we
   identify $f(F_A)$ with $\alpha$.

   It is a standard result of Chern-Weil theory that $f(F_A)\in \Omega^{2k}(M;\mathbb{R})$ is closed, that is,
   $df(F_A)=0$ and that the class of $f(F_A)$ in $ H^*_{dR}(M)$  does not depend on $A$,
   that is,
   if $A'$ is another connection, $f(F_A)-f(F_{A'})$ is exact.
   In particular, let $a\in \Omega^1(P;\mathcal{G})$ be horizontal and equivariant.
   Then $A+a$ is another connection and we have(cf.\ p.297, \cite{Nomizu})
   $$f(F_{A+a})-f(F_A)=k\ d\left(\int_0^1f(a,F_{A+ta},\cdots, F_{A+ta})dt\right).$$
   This implies that there is an algebra homomorphism $I(\mathcal{G})\rightarrow H^{2*}_{dR}(M)$
   determined by the principal bundle $P$.

   \section{The Maurer-Cartan 3-form of a connected compact Lie group}

We begin by the identification of the de Rham cohomology of a compact connected Lie group
with the algebra consisting of the {\it bi-invariant} forms.
We will also explain the notational conventions in the expressions
such as $\langle (A-A_0)\wedge [(A-A_0)\wedge(A-A_0)]\rangle$
which appeared repeatedly in the introduction and will continue to do so
throughout the paper.

 Let $G$ be a connected compact Lie group and
$\mathcal{H}^*(G)$ denote the set of all {\it bi-invariant} real valued forms
on $G$. That a form $\alpha\in \Omega^*(G;\mathbb{R})$ is bi-invariant means
that it satisfies $L_g^*\alpha=R_g^*\alpha=\alpha$ for any $g\in G$. Bi-invariant forms are
closed and
we have(Theorem 12.1, \cite{Chevalley}):
\bdm \mathcal{H}^*(G)\equiv H^*_{\rm dR}(G).\tag{3-1}\edm

Let $q:\tilde G\rightarrow G$ be a finite cover by a connected group $\tilde G$.
Since $q$ is in fact a group homomorphism, we may easily prove that,
for any  bi-invariant $\alpha\in \Omega^*(G;\mathbb{R})$,
$q^*\alpha$ is also bi-invariant.
It is also straightforward to see
that any bi-invariant form of $\tilde G$ is the pull-back of a unique bi-invariant
form of $G$. Therefore (3-1) implies the isomorphism
\bdm H^*_{\rm dR}(G)\equiv H^*_{\rm dR}(\tilde G).\tag{3-2}\edm

It is well-known that $G\cong \tilde G /Z$ where $\tilde G$ is the product group
$$T^k\times G_1\times \cdots \times G_l$$
of the torus group $T^k$ and simply connected simple Lie groups $G_i$, $i=1,\cdots, l$
and $Z$ is a finite subgroup
of the center of $\tilde G$. Therefore considering (3-2) we have:
$$ H^*_{\rm dR}(G)\cong H^*_{\rm dR}(T^k)\otimes H^*_{\rm dR}(G_1)\otimes\cdots\otimes H^*_{\rm dR}(G_l) .$$

Recall the Maurer-Cartan form $\theta\in \Omega^1(G;\mathcal{G})$ and the equalities
\bdm L_g^*\theta=\theta, \ R_g^*\theta=Ad_{g^{-1}}\circ \theta \tag{3-3}\edm
for any $g\in G$.
We consider the Maurer-Cartan 3-form on $G$
\bdm \Theta=-\frac{1}{6}\langle\theta\wedge[\theta\wedge\theta]\rangle \tag{3-4}\edm
where $\langle\cdot,\cdot\rangle$ is an Ad-invariant inner product on $\mathcal{G}$.
Here we understand $\theta\wedge\theta$ as a $\mathcal{G}\otimes \mathcal{G}$-valued
$2$-form defined by
$$(\theta\wedge\theta)(X, Y)=\theta(X)\otimes\theta(Y)-\theta(Y)\otimes\theta(X)$$
for any $X, Y \in T_gG$ and any $g\in G$.
Furthermore we identify any bilinear function on $\mathcal{G}\times\mathcal{G}$ with a linear
map on the vector space $\mathcal{G}\otimes\mathcal{G}$. With these conventions being understood,
any expression such as the RHS of (3-4) makes a good sense.
By exploiting (3-3) it is rather straightforward to see that $\Theta$ is a bi-invariant $3$-form and in particular
that $\Theta$ is closed.

Consider the case $G=S^3\subset \mathbb{H}$.
Then $\mathcal{G}=T_1S^3$ has the standard basis $\{i, j, k\}$.
Let us consider $\mathcal{G}\subset T_1\mathbb{H}=T_{(1,0,0,0)}\mathbb{R}^4$ with the Euclidean metric.
Then we have:
$$\Theta(i,j,k)=-\langle i, [j, k]\rangle=-2.$$
Thus $\Theta$ is not zero and $[\Theta]$ is a basis for $H^3_{\rm dR}(S^3)$.

Let $\varphi: H\rightarrow K$ be a homomorphism between Lie groups.
If $\theta_H$ and $\theta_K$ are the Maurer-Cartan forms,
we have:
$$\varphi^*\theta_K=d\varphi_e\circ \theta_H.$$
Furthermore assume $d\varphi_e:\mathcal{H}\rightarrow\mathcal{K}$ is an injective
homomorphism between the Lie algebras.
Assume $\mathcal{H}$
is equipped with the inner product which is the pull-back of an Ad-invariant inner product on $\mathcal{K}$.
Then the inner product on $\mathcal{H}$ is also Ad-invariant. Now let
$\Theta_H$ and $\Theta_K$ be the Maurer-Cartan 3-forms of $H$ and $K$ respectively.
Then we have:
\bdm \varphi^*\Theta_K=\Theta_H .\tag{3-5}\edm

Let $G$ be a compact simply connected simple Lie group. It is well known that there is
a homomorphism
$$\varphi:{\rm SU}(2)\rightarrow G$$
such that the homomorphism $d\varphi_e: su(2)\rightarrow \mathcal{G}$ between
the Lie algebras is injective. Note that there is a unique Ad-invariant inner product on
$\mathcal{G}$ up to multiplication by a positive real number.
It is also a classical fact(see for instance \S 1, \cite{Kachi}) that
\bdm
H^2(G;\mathbb{R})= 0  {\rm\ \ and\ \ } H^3(G;\mathbb{R})\cong \mathbb{R}. \tag{3-6}
\edm
In particular equation (3-5), together with the fact
the Maurer-Cartan 3-form is a basis for $H^3_{\rm dR}({\rm SU}(2))\cong H^3_{\rm dR}(S^3)$,
implies that the class $[\Theta]$
of the Maurer-Cartan $3$-form is a basis
for $H^3_{\rm dR}(G)$.

Now assume
$G= G_1\times\cdots \times G_l$ where $G_k$, $k=1,\cdots, l$, are simply connected simple Lie groups.
Then
we have the decomposition
$$\mathcal{G}\equiv \mathcal{G}_1\oplus\cdots\oplus\mathcal{G}_l$$
for the Lie algebras. An Ad-invariant inner product on $\mathcal{G}$\
makes this decomposition orthogonal and therefore
it is determined by its restrictions to $\mathcal{G}_k$'s.
Thus without loss of generality we may assume each $\mathcal{G}_k$ is equipped with an Ad-invariant inner product
and $\mathcal{G}$ is given the direct sum inner product. Use these inner products
to define the Maurer-Cartan $3$-forms $\Theta$ and $\Theta_k$'s respectively of $G$ and $G_k$'s.
Then we have
\bdm \Theta=q_1^*\Theta_1 +\cdots +q_l^*\Theta_l \tag{3-7}\edm
where $q_k:G\rightarrow G_k$, $k=1,\cdots, l$, are the projections.
On the other hand since $H^i_{\rm dR}(G_k)=0$ for $i=1,2$ we have that
$$H^3_{\rm dR}(G)\equiv  H^3_{\rm dR}(G_1)\oplus\cdots\oplus H^3_{\rm dR}(G_l).$$
Note that there is a unique Ad-invariant inner product on each
$\mathcal{G}_k$ up to multiplication by a positive real number.
Therefore (3-7) implies that the set of all $[\Theta]$, each of which is determined by a choice
of an Ad-invariant inner product
on $\mathcal{G}$, is an open cone in $H^3_{\rm dR}(G)$.

Slightly more generally let $G$ be a compact connected semi-simple Lie group.
Then $G\cong \tilde G/Z$ where $\tilde G= G_1\times\cdots G_l$ for some simply connected simple Lie groups $G_k$,
$k=1,\cdots, l$,
and $Z$ is a subgroup of the center of $\tilde G$.
 By applying (3-2) the set of all classes $[\Theta]$ of Maurer-Cartan $3$-forms
is still an open cone in $H^3_{\rm dR}(G)$.
Since the image of the homomorphism
 $H^3(G; \mathbb{Z})\rightarrow H^3(G;\mathbb{R})\equiv H^3_{\rm dR}(G)$
is a full lattice, we conclude that we may choose an Ad-invariant metric on $\mathcal{G}$ so that
$[\Theta]\in H^3_{\rm dR}(G)$ is integral.

Now let $G$ be any compact connected Lie group.
Then $G'=G/Z_0$ is semi-simple(or trivial),
where $Z_0$ is the identity component of
the center of $G$. We have a natural decomposition of the Lie algebra,
$\mathcal{G}=\mathcal{Z}\oplus \mathcal{G}'$, where
$\mathcal{Z}$ is the center of $\mathcal{G}$ and
$\mathcal{G}'=[\mathcal{G},\mathcal{G}]$.
Let $q:G\rightarrow G'$ be the quotient homomorphism. Any Ad-invariant inner product
on the Lie algebra of $G'$ is given,
we may equip $\mathcal{G}$ with an
Ad-invariant metric so that the restriction of $dq_e$ to $\mathcal{G}'$
is an isometry. Now let $\Theta$ and $\Theta'$ be the Maurer-Cartan $3$-forms of
$G$ and $G'$ respectively. Then we have:
\bdm \Theta= q^*\Theta'.\tag{3-8}\edm
To see this equality, note that $q^*\Theta'$ is bi-invariant. Therefore it is
enough to see that the equality holds when evaluated
at a triple $(X,Y,Z)\in \mathcal{G}^3$. Thus it is sufficient to
see that $\langle X,[Y,Z]\rangle=\langle dq_eX,[dq_eY,dq_eZ]\rangle$ holds for any
$(X,Y,Z)\in \mathcal{G}^3$, which is straightforward.

In particular we have proved the following fact:
\begin{theorem}
Let $G$ be a compact connected Lie group. Then there is an Ad-invariant inner product
on the Lie algebra of $G$ so that the cohomology class of the Maurer-Cartan $3$-form is
integral.
\end{theorem}

\section{The Chern-Simons functional defined by a section}
In case of trivializable bundle the Chern-Simons functional is usually given by
choosing a section. Here we observe that the Chern-Simons functional
defined by choosing a reference connection is indeed a generalization
of this special case.
Furthermore we prove (1-2), (1-3) and (1-4) for this case, which is somewhat easier
and can be regarded as a guidance to the full proofs for the more general case.

It is well known that all Chern-Weil forms are exact as forms on the principal
bundle(\cite{Simons}). Let us keep the notations of the introduction.
Then in particular we have:

\begin{lemma}
Let $A\in \mathcal{A}(P)$ and $\langle \cdot, \cdot\rangle$ be an
Ad-invariant inner product on $\mathcal{G}$. Then we have that
$$d(\langle A \wedge F_A\rangle -\frac{1}{6}\langle A\wedge[A\wedge A]\rangle)=\langle F_A \wedge F_A \rangle.$$
\end{lemma}
\bp
Apply the identities
$$dA=F_A -\frac{1}{2}[A\wedge A], \ dF_A = -[A\wedge F_A]$$
to conclude
$$d\langle A\wedge F_A\rangle=\langle F_A \wedge F_A\rangle +\frac{1}{2}\langle F_A \wedge [A\wedge A] \rangle .$$
In fact in the above we also used the identity
$$\langle \alpha \wedge [\beta\wedge \gamma] \rangle
=(-1)^{r(p+q)}\langle \gamma \wedge [\alpha \wedge \beta ]\rangle$$
if $\alpha, \beta, \gamma \in \Omega^*(P; \mathcal{G})$ are homogeneous
forms of degree $p, q, r$ respectively, which is a consequence
of Ad-invariance of the inner product $\langle\cdot,\cdot\rangle$.

We also have
$$d\langle A \wedge [A\wedge A]\rangle =3\langle dA \wedge [A\wedge A]\rangle.$$
Use the identity $dA=F_A -\frac{1}{2}[A\wedge A]$ again and observe
$$\langle [A\wedge A] \wedge [A\wedge A]\rangle = \langle A\wedge [A \wedge [A\wedge A]]\rangle=0 .$$
In fact the equality $[A \wedge [A\wedge A]]=0$ holds, which is a consequence of the Jacobi identity.
Therefore we conclude that
$$d\langle A \wedge [A\wedge A]\rangle=3\langle F_A \wedge [A\wedge A]\rangle , $$
which immediately proves the lemma.
\ep

Now we assume that the dimension of $M$ is $3$ and $M$ is oriented and closed.
And we assume furthermore that the principal bundle
$\pi:P \rightarrow M$ is trivializable and $s:M\rightarrow P$ is a section.
Write $\alpha(A)=\langle A \wedge F_A\rangle -\frac{1}{6}\langle A\wedge[A\wedge A]\rangle$.
Then we define a Chern-Simons functional $CS_s:\mathcal{A}(P)\rightarrow \R $ by
$$CS_s(A)=\int_M s^*\alpha(A)$$
for any $A\in\mathcal{A}(P)$. Note that $\alpha(A)$ need not be horizontal
even if it is invariant.
Therefore pulling it back by $s:M\rightarrow P$ is a necessity.

On the other hand note that the section $s$ defines a unique flat connection $A_0\in \mathcal{A}(P)$
such that $s(M)$ is an integral submanifold of the horizontal distribution given by $A_0$.
Let $A\in \mathcal{A}(P)$. Then we consider the $3$-form:
$$cs(A)=\langle F_A\wedge (A-A_0)\rangle - {\frac{1}{6}} \langle (A-A_0)\wedge [(A-A_0)\wedge(A-A_0)]\rangle.$$
Then we have that
$$s^*\alpha(A)=cs(A),$$
where we regard $cs(A)$ as a form on $M$. Thus the functional given by a global section is a special
case of the functional determined by choosing a reference connection.

Now we observe the following:

\begin{lemma}
Let $A\in\mathcal{A}(P)$. Then we have that
$${dCS_s}_A(a)=2\int_M\langle F_A\wedge a\rangle$$
for any $a\in \bar\Omega^1(P;\mathcal{G})$.
\end{lemma}
\bp
Note that $F_{A+ta}=F_A+td_Aa +\frac{1}{2}t^2[a\wedge a]$.
Therefore we have:
$$\left.\frac{d}{dt}\right|_{t=0}\alpha(A+ta)=\langle a\wedge F_A\rangle
+\langle A\wedge d_Aa\rangle -\frac{1}{2}\langle a\wedge[A\wedge A]\rangle.$$
Since $d_Aa=da +[A\wedge a]$ and $F_A=dA+\frac{1}{2}[A\wedge A]$,
we may rewrite this as follows:
$$
\aligned
\left.\frac{d}{dt}\right|_{t=0}\alpha(A+ta)&=\langle a\wedge F_A\rangle +\langle da\wedge A\rangle
+\frac{1}{2}\langle a\wedge [A\wedge A]\rangle\\
&=2\langle a\wedge F_A\rangle + d\langle a\wedge A\rangle.
\endaligned
$$

Note that $\langle a\wedge F_A\rangle$ is horizontal and invariant.
Therefore we have:
$$
\aligned
\left.\frac{d}{dt}\right|_{t=0}CS_s(A+ta)&=2\int_M s^* \langle a\wedge F_A\rangle\\
&=2\int_M \langle a\wedge F_A\rangle,
\endaligned
$$
as desired.
\ep

The RHS of the equation in the lemma above does not depend on the section $s$.
Therefore, if $s':M\rightarrow P$ is another section, $CS_{s'}-CS_s:\mathcal{A}(P)\rightarrow\mathbb{R}$
is a constant.
Let $\varphi:P\rightarrow P$ be a gauge transformation. Then the equality,
$CS_s(\varphi^*A)=CS_{\varphi\circ s}(A)$, holds. Therefore
$CS_s(\varphi^*A)-CS_s(A)$ must be constant of $A$.

The following leads to a concrete value of the constant $CS_s(\varphi^*A)-CS_s(A)$ and in fact
implies that it depends only on the path component of $\varphi$ in the gauge group ${\rm GA}(P)$.
\begin{lemma}
Let $u: P\rightarrow G$
be the map defined by the equation, $\varphi(p)=pu(p)$, for any $p\in P$ and
let $\theta\in \Omega^1(G;\mathcal{G})$ and $\Theta\in \Omega^3(G;\mathbb{R})$
denote respectively the Maurer-Cartan form and 3-form.  Then we have:
$$\alpha(\varphi^*A)-\alpha(A)=u^*\Theta
                  +d\langle Ad_{u^{-1}}A\wedge u^*\theta\rangle.$$
\end{lemma}
\bp
Note that $\varphi^*A=Ad_{u^{-1}}A +u^*\theta$ and $F_{\varphi^*A}=Ad_{u^{-1}}F_A$.
Therefore we have:
$$\aligned
&\alpha(\varphi^*A)- \alpha(A) + \frac{1}{6}\langle u^*\theta \wedge[u^*\theta \wedge u^*\theta]\rangle\\
=&\langle u^*\theta \wedge Ad_{u^{-1}}F_A\rangle
 -\frac{1}{2}\langle Ad_{u^{-1}}A\wedge [u^*\theta \wedge u^*\theta]\rangle
-\frac{1}{2}\langle u^*\theta\wedge[Ad_{u^{-1}}A \wedge Ad_{u^{-1}}A]\rangle .
 \endaligned
$$

On the other hand let us begin with applying the Leibniz rule:
$$ d\langle Ad_{u^{-1}}A\wedge  u^*\theta\rangle=\langle d (Ad_{u^{-1}}A)\wedge u^*\theta\rangle
                 -\langle Ad_{u^{-1}}A\wedge d u^*\theta\rangle.$$

Firstly note that
$du^*\theta=u^*d\theta=-\frac{1}{2}u^*[\theta\wedge\theta]=-\frac{1}{2}[u^*\theta\wedge u^*\theta]$.

Secondly we have
$$\aligned
d(Ad_{u^{-1}}A)&=(d Ad_{u^{-1}})\wedge A + Ad_{u^{-1}}dA\\
     &=(Ad_{u^{-1}} ad({u^{-1}}^*\theta))\wedge A +Ad_{u^{-1}}(F_A -\frac{1}{2}[A\wedge A])\\
     &=Ad_{u^{-1}}[{u^{-1}}^*\theta\wedge A] +Ad_{u^{-1}}F_A -\frac{1}{2}[Ad_{u^{-1}} A\wedge Ad_{u^{-1}} A] .
\endaligned
$$
Here $dAd_{u^{-1}}$ should be understood as the derivative of a vector valued function,
$Ad_{u^{-1}}:P\rightarrow End(\mathcal{G})$.
Also note the equality ${u^{-1}}^*\theta=-Ad_u u^*\theta$. Then we have:
$$\aligned
&\langle d(Ad_{u^{-1}}A)\wedge u^* \theta\rangle\\
=&-\langle (Ad_{u^{-1}}[Ad_u u^*\theta \wedge A])\wedge u^*\theta\rangle
   +\langle Ad_{u^{-1}}F_A \wedge u^*\theta\rangle\\
 &\ \ \ \ \ \ \ \ \ \ \ \ \ \ \ \  -\frac{1}{2}\langle[Ad_{u^{-1}}A\wedge Ad_{u^{-1}}A]\wedge u^*\theta\rangle\\
=&-\langle Ad_{u^{-1}}A\wedge [ u^*\theta \wedge u^*\theta]\rangle
 +\langle Ad_{u^{-1}}F_A\wedge u^*\theta\rangle
-\frac{1}{2}\langle u^*\theta\wedge [Ad_{u^{-1}}A\wedge Ad_{u^{-1}}A] \rangle.
\endaligned
$$

Thus we obtain the equality:
$$
\aligned
& d\langle Ad_{u^{-1}}A\wedge  u^*\theta\rangle\\
=& \langle u^*\theta \wedge Ad_{u^{-1}}F_A\rangle
 -\frac{1}{2}\langle Ad_{u^{-1}}A\wedge [u^*\theta \wedge u^*\theta]\rangle
 -\frac{1}{2}\langle u^*\theta\wedge[Ad_{u^{-1}}A \wedge Ad_{u^{-1}}A]\rangle,
\endaligned
$$
as desired.
\ep

Let $s:M\rightarrow P$ be a section and let $\varphi:P\rightarrow P$ and
$u:P\rightarrow G$ be as in the lemma above. We define $\deg \varphi$ by
$$\deg \varphi = \int_M (u\circ s)^*\Theta \in \mathbb{R}.$$
By 3.1 above, if $G$ is connected and compact, we may assume the Ad-invariant metric
on $\mathcal{G}$
is such that $[\Theta]\in H^3_{\rm dR}(G)$ is integral. Then we have that
$$\deg \varphi = \langle (u\circ s)^*[\Theta], [M]\rangle \in \mathbb{Z}.$$

We need to observe:
\begin{lemma}
If $G$ is connected and compact, the invariant $\deg \varphi$ does not depend
on the section $s:M\rightarrow P$.
\end{lemma}
\bp
Let $s':M\rightarrow P$ be another section. Then there is a smooth map $w:M\rightarrow G$
such that $s'(x)=s(x)w(x)$ for any $x\in M$. And we have
$(u\circ s')(x)=u(s'(x))=u(s(x)w(x))=w(x)^{-1}u(s(x))w(x)$
for any $x\in M$.
Now from the claim below
we have $(u\circ s')^*[\Theta]=(u\circ s)^*[\Theta]$ and the assertion of the lemma follows.

{\it Claim.} Let $X$ be a smooth manifold and
let $f, g: X\to G$ be smooth maps and $h:X\to G$ be defined by $h(x)=f(x)g(x)$ for any $x\in X$.
Then we have $h^*[\Theta]=f^*[\Theta] +g^*[\Theta] \in H^3_{dR}(X)$.

{\it Proof.} Note that $h$ is the composite
$$G\overset{\Delta}{\rightarrow} G\times G\overset{f\times g}{\longrightarrow} G\times G\overset{\mu}{\rightarrow}G,$$
where $\Delta$ is the diagonal and $\mu$ is the multiplication.

Firstly assume
that $G$ is semisimple. Then we have that $H^3_{\rm dR}(G\times G)\equiv H^3_{\rm dR}(G)\otimes \mathbb{R} \oplus \mathbb{R}\otimes H^3_{\rm dR}(G)$
since $H^i_{\rm dR}(G)=0$ if $i=1,2$, and $H^0_{\rm dR}(G)=\mathbb{R}$. It follows that $\mu^*a=a\times 1+1\times a$ for any $a\in H^3_{\rm dR}(G)$.
Therefore we have that $h^*a=f^*a +g^*a$ for any $a\in H^3_{\rm dR}(G)$. For a general compact connected group $G$,
consider the projection $q:G\rightarrow G'$ where $G'=G/Z_0$ is the quotient by the identity component of
the center. Then $qh=(qf)(qg)$ and we have $(qh)^*[\Theta']=(qf)^*[\Theta']+(qg)^*[\Theta']$.
Now by (3-8) we have the equality.
\ep

Thus $\deg \varphi\in \mathbb{R}$ is well-defined if $P\rightarrow M$ is trivializable
and the structure group $G$ is connected and compact.
By applying 4.3 above we conclude the well-known fact:
\begin{theorem}
Let $P\rightarrow M$ be a trivializable bundle with a connected compact structure group
over an oriented closed $3$-manifold $M$.
Then for any section $s:M\rightarrow P$, any $A\in\mathcal{A}(P)$ and
any gauge transformation $\varphi:P\rightarrow P$, we have:
$$CS_s(\varphi^*A)-CS_s(A)=\deg \varphi .$$
\end{theorem}

\section{The Chern-Simons functional defined by a reference connection} 
This section begins with justifications of (1-2) and (1-3).
Then we provide a detailed review of the properties of the $1$-form $dCS$
on $\mathcal A (P)$
which leads to
an understanding of the general behavior of the Chern-Simons functional
under gauge transformations.

Let $\pi:P\rightarrow M$ be a principal $G$-bundle over a manifold $M$.
Let $A_0, A_1\in \mathcal{A}(P)$. Then from the Chern-Weil theory(see \S 2 above) we have
$$
\aligned
\langle F_{A_1}\wedge & F_{A_1}\rangle - \langle F_{A_0}\wedge F_{A_0}\rangle\\
 &=d(2\int_0^1\langle(A_1 -A_0)\wedge F_{A_0 +t(A_1 -A_0)}\rangle dt).
\endaligned
$$
We fix an $A_0\in \mathcal{A}(P)$ and define for any $A\in \mathcal{A}(P)$
the 3-form
$$cs(A)=2\int_0^1\langle(A -A_0)\wedge F_{A_0 +t(A -A_0)}\rangle dt.$$

Recall that for any $a\in \bar\Omega^1(P;\mathcal{G})$ we have: $F_{A+a}=F_A +d_Aa+\frac{1}{2}[a\wedge a]$.
Then we have other expressions for $cs(A)$ as well:
$$
\aligned
cs(A)&=2\langle(A-A_0)\wedge  F_{A_0}\rangle +\langle(A-A_0)\wedge d_{A_0}(A-A_0)\rangle\\
     &\ \ \ \ \ \ \ \ \ \ \ \ \ \ \ \ +\frac{1}{3}\langle(A-A_0)\wedge[(A-A_0)\wedge(A-A_0)]\rangle\\
     &=\langle (A-A_0)\wedge(F_A+F_{A_0})\rangle-\frac{1}{6}\langle (A-A_0)\wedge[(A-A_0)\wedge(A-A_0)]\rangle.
\endaligned
$$
This shows (1-2) holds.

On the other hand we have
$$
\aligned
\left.\frac{d}{dt}\right|_{t=0}cs(A+ta)=\langle a\wedge(F_A +&F_{A_0})\rangle +\langle(A-A_0)\wedge d_A a\rangle\\
                                                 &-\frac{1}{2}\langle a\wedge[(A-A_0)\wedge(A-A_0)]\rangle.
\endaligned
$$
Using the identities $d_Aa=d_{A_0}a+[(A-A_0)\wedge a]$ and
$F_A=F_{A_0}+d_{A_0}(A-A_0)+\frac{1}{2}[(A-A_0)\wedge (A-A_0)]$ subsequently,
we have:
$$
\aligned
\left.\frac{d}{dt}\right|_{t=0}cs(A+ta)&=\langle a\wedge F_A\rangle +\langle a\wedge F_{A_0}\rangle +\langle(A-A_0))\wedge d_{A_0}a\rangle\\
                           &\ \ \ \ \ \ \ \ \ \ \ \ \ +\frac{1}{2}\langle a\wedge[(A-A_0)\wedge(A-A_0)]\rangle\\
                           &=2\langle a\wedge F_A\rangle +\langle d_{A_0}a\wedge(A-A_0)\rangle -\langle a\wedge d_{A_0}(A-A_0)\rangle\\
                           &=2\langle a\wedge F_A\rangle +d\langle a\wedge(A-A_0)\rangle.
\endaligned
$$
This immediately proves the equality, $d CS_A(a)=2\int_M\langle F_A\wedge a\rangle$, which is (1-3).

Note that the group ${\rm GA}(P)$ of all gauge transformations of $P$ acts on $\mathcal{A}(P)$.
The following has been asserted in \cite{Dostoglou}, where a justification has not been provided.
We add a detailed proof.

\begin{lemma}
The 1-form $dCS$ on $\mathcal{A}(P)$ is invariant and horizontal. That is, we have:
$$\left.\frac{d}{dt}\right|_{t=0}CS(\varphi^*(A+ta))=\left.\frac{d}{dt}\right|_{t=0}CS(A+ta)
\ \ {\rm and}\ \ \left.\frac{d}{dt}\right|_{t=0}CS(\varphi_t^*A)=0,$$
for any $A\in \mathcal{A}(P)$, any $a\in \bar\Omega^1(P;\mathcal{G})$, any gauge transformation
$\varphi$ and any smooth 1-parameter family $\varphi_t$, $-\epsilon <t <\epsilon$, of gauge
transformations such that $\varphi_0=1$.
\end{lemma}
\bp
Let $u:P\rightarrow G$ be the map determined by
$\varphi$. Since $a$ is equivariant, we have $\varphi^*a=Ad_{u^{-1}} a$. It follows that
$$
\aligned
\left.\frac{d}{dt}\right|_{t=0}CS&(\varphi^*(A+ta))=\left.\frac{d}{dt}\right|_{t=0}CS(\varphi^*A+t\varphi^*a)\\
                                 &=2\int_M\langle F_{\varphi^*A}\wedge \varphi^*a\rangle\\
                                 &=2\int_M\langle (Ad_{u^{-1}} F_A)\wedge (Ad_{u^{-1}} a)\rangle\\
                                 &=2\int_M \langle F_A\wedge a\rangle.
\endaligned
$$
This proves the first equality.

Let $u_t$ denote the map $P\rightarrow G$ determined by $\varphi_t$ for each $t$.
We have that $\varphi_t^*A=A+u_t^{-1}d_Au_t$, where $u^{-1}d_Au$ means $(u^{-1}du)^H$
if $u:P\rightarrow G$ is any map and $H$ is the horizontal distribution determined by $A$.
Let $X:P\rightarrow \mathcal{G}$ denote the map
defined by $X(p)=(u_t(p))^{\centerdot}_{t=0}$. We observe:

{\it Claim:} $\left.\frac{d}{dt}\right|_{t=0}u_t^{-1}du_t =dX$.

{\it Proof.} Let $Y_p\in T_p P$ and let $\gamma:(-\delta,\delta)\rightarrow P$ be a curve
such that $\dot{\gamma}(0)=Y_p$. Then define
$f:(-\epsilon, \epsilon)\times(-\delta,\delta)\rightarrow G$ by $f(t,s)=u_t(p)^{-1}u_t(\gamma(s))$.
We may write:
$$(\left.\frac{d}{dt}\right|_{t=0}u_t^{-1}du_t)(Y_p)
    =\left.\frac{d}{dt}\right|_{t=0}df(\left.\frac{\partial}{\partial s}\right|_{(t,0)})\in T_eG=\mathcal{G}.$$
In general, let $g:(-\epsilon, \epsilon)\times(-\delta,\delta)\rightarrow N$ be any $C^2$-map into a smooth
manifold such that
there is an $x\in N$ for which $g(t,0)=x$, $g(0,s)=x$ holds for any $t$, $s$. Then we have:
$$\left.\frac{d}{dt}\right|_{t=0}dg(\left.\frac{\partial}{\partial s}\right|_{(t,0)})
 =\left.\frac{d}{ds}\right|_{s=0}dg(\left.\frac{\partial}{\partial t}\right|_{(0,s)})\in T_xN.$$
Therefore the claim is proved by the following:
$$
\aligned
\left.\frac{d}{dt}\right|_{t=0}df(\left.\frac{\partial}{\partial s}\right|_{(t,0)})
&=\left.\frac{d}{ds}\right|_{s=0}df(\left.\frac{\partial}{\partial t}\right|_{(0,s)})\\
&=\left.\frac{d}{ds}\right|_{s=0}(-X(p) +X(\gamma(s)))=dX(Y_p).
\endaligned
$$

Now note that $d_AF_A=0$. Then the following proves the second equality:
$$
\aligned
\left.\frac{d}{dt}\right|_{t=0}CS(\varphi_t^*A)&=\left.\frac{d}{dt}\right|_{t=0}CS(A+ u_t^{-1}d_Au_t)\\
                                                &=2\int_M\langle F_A\wedge d_A X\rangle =2\int_M d\langle F_A\otimes X\rangle=0.
\endaligned
$$
\ep

The first equality of the lemma above implies that $CS(\varphi^*A)-CS(A)$ does not depend on $A$.
The second means in fact that $\frac{d}{dt}(CS(\varphi_t^*A))=0$ for any smooth
$1$-parameter family of gauge transformations $\varphi_t$ and therefore it follows that $CS(\varphi^*A)-CS(A)$ is constant
if $\varphi$ varies within a path component of ${\rm GA}(P)$.
Since the choice of the reference connection $A_0$ affects the Chern-Simons functional only by addition
of a constant, $CS(\varphi^*A)-CS(A)$ does not depend
on the reference connection either.  Therefore we have:
\begin{proposition}
The real number, $CS(\varphi^*A)-CS(A)$, depends only on the path component of $\varphi$ in ${\rm GA}(P)$.
In particular
it depends neither on $A$ nor on the reference connection.
\end{proposition}

\section{Cohomology of flat bundles} 
In this section we describe the cohomology rings of a flat bundle
and its adjoint bundle of Lie groups using the cohomology rings of the base manifold
and the structure group $G$ under the assumption that $G$ is connected and compact.

Let $\pi:P\rightarrow M$ be a principal bundle whose structure group $G$ is connected and compact .
Let $A\in \mathcal{A}(P)$ and $H$ be the horizontal distribution
of $P$ determined by $A$. Let $V$ denote the vertical distribution of $P$.
Then we let
$$\pi_A: TP\rightarrow V$$
denote the projection given by
the decomposition $TP=V\oplus H$.

For any $p\in P$  we have a map $\kappa_p:P_x\rightarrow G$, $x=\pi(p)$,
defined by $\kappa_p(pg)=g$ for any $g\in G$. Let $h\in G$. Then we have $\kappa_{ph}=L_{h^{-1}}\kappa_p$, which follows
from the identities:
$\kappa_{ph}(pg)=\kappa_{ph}(ph(h^{-1}g))=h^{-1}g=L_{h^{-1}}\kappa_p(pg)$.
Recall the set $\mathcal{H}^*(G)$ of all bi-invariant real valued forms. Then the equality,
$\kappa_{ph}=L_{h^{-1}}\kappa_p$ implies that, for any
$\alpha\in \mathcal{H}^*(G)$, the pull-back $\kappa_p^*\alpha \in \Omega^*(P_x;\mathbb{R})$
does not depend on $p\in P_x$.

Now we define $\hat\alpha\in \Omega^k(P; \mathbb{R})$, for any $\alpha\in \mathcal{H}^k(G)$,
by the rule:
\bdm
\hat\alpha(X_1, \cdots X_k)=(\kappa_p^*\alpha)(\pi_A X_1, \cdots , \pi_A X_k),\tag{6-1}
\edm
for any tangent vectors $X_1,\cdots, X_k\in T_{p'} P$, $p'\in P$, choosing any $p\in P_{\pi(p')}$.
We write
$$E_A: \mathcal{H}^*(G)\rightarrow \Omega^*(P;\mathbb{R})$$
for the map
defined by $E_A(\alpha)=\hat\alpha$ for any $\alpha\in\mathcal{H}^*(G)$.

Now we observe:
\begin{lemma}
If $A$ is flat, $E_A(\alpha)\in \Omega^*(P; \mathbb{R})$ is closed for any $\alpha\in \mathcal{H}^*(G)$.
\end{lemma}
\bp
Let $U\subset M$ be an open set for which there is a section $s:U\rightarrow P$
horizontal with respect to $A$.

Define $\kappa_s :P_U=\pi^{-1}U\rightarrow G$ by $p=s(\pi(p))\kappa_s(p)$ for any $p\in P_U$. Then we have
$$\hat \alpha|_{P_U}=\kappa_s^*\alpha.$$
In fact this equality follows from the identity
$$d\kappa_{s(x)}\pi_A=d\kappa_s :T_pP\rightarrow T_{\kappa_s(p)}G,$$
for any $p\in P$ where $x=\pi(p)$. This is easily seen if one notes any $Y\in T_pP$ is the velocity
at $t=0$ of a curve $s(\delta(t))\gamma(t)$ on $P$ where $\delta$ and $\gamma$ are some curves respectively on $M$ and on $G$
such that $\delta(0)=x$ and $\gamma(0)=\kappa_{s(x)}(p)$. Then both $d\kappa_{s(x)}\pi_A(Y)$ and $d\kappa_s(Y)$
are $\dot\gamma(0)$.

Therefore we conclude $d\hat\alpha = 0$ on $P_U$. Since $A$ is flat, $M$ is covered by
open sets on which horizontal sections exist. This proves that $\hat \alpha =E_A(\alpha)$ is closed
on $P$.
\ep

Thus assuming $A$ is flat we have a map $H^*_{\rm dR}(G)\rightarrow H^*_{\rm dR}(P)$
induced by $E_A$ since $\mathcal{H}^*(G)\equiv H^*_{\rm dR}(G)$, which we still denote
by $E_A$. It is clear that $E_A:H^*_{\rm dR}(G)\rightarrow H^*_{\rm dR}(P)$ is a homomorphism
between algebras.

Consider the inclusion $\iota_p :G \overset{\kappa_p^{-1}}{\longrightarrow} P_{\pi(p)}\subset P$
for any $p\in P$. In fact $\iota_p:G\rightarrow P$ is given by $\iota_p(g)=pg$, $g\in G$. Then
we have that $\iota_p^*\hat \alpha=\alpha$ for any $p\in P$. Therefore we have proved the following:
\begin{theorem}
Let $P\rightarrow M$ be a flat bundle over a smooth manifold with a compact connected
structure group $G$. Then for any flat connection $A$ on $P$ there is a homomorphism between
algebras,
$$E_A:H^*_{\rm dR}(G)\rightarrow H^*_{\rm dR}(P),$$
such that, for any $p\in P$, $\iota_p^*E_A$
is the identity on $H^*_{\rm dR}(G)$.
\end{theorem}

The above means a flat bundle with a compact connected structure group $G$ is a fiber bundle
to which the Leray-Hirsch theorem(cf.\ \cite{Hatcher}) applies: We have an isomorphism of $H^*(M)$-modules
$$H^*_{\rm dR}(P)\cong H^*_{\rm dR}(M)\otimes H^*_{\rm dR}(G),$$
where the isomorphism $H^*_{\rm dR}(M)\otimes H^*_{\rm dR}(G)\rightarrow H^*_{\rm dR}(P)$  is given by
$a\otimes b\rightarrow (\pi^*a)\wedge (E_A(b))$.
Moreover we consider $H^*_{\rm dR}(M)\otimes H^*_{\rm dR}(G)$ with the ring structure
given by
$$(a_1\otimes b_1)(a_2\otimes b_2)=(-1)^{kl}(a_1\wedge a_2)\otimes(b_1\wedge b_2),$$
if $a_1, a_2\in H^*_{\rm dR}(M)$ and $b_1, b_2\in H^*_{\rm dR}(G)$ and in particular
$b_1$ and $a_2$ are homogeneous classes of degrees $k$ and $l$ respectively.
Then since $E_A$ is a homomorphism of algebras,
we have:
\begin{corollary}
Let $P\rightarrow M$ and $G$ be as in 6.2 above. Then
there is an isomorphism of algebras:
$$H^*_{\rm dR}(P)\cong H^*_{\rm dR}(M)\otimes H^*_{\rm dR}(G).$$
\end{corollary}

We may proceed with the adjoint bundle $ Ad\ P $ of Lie groups in a similar way
as can be seen in the below.

Let $A$ be any connection
on $P$ and $H$ be the horizontal distribution determined by $A$. Let
$q:P\times G \rightarrow Ad\ P =(P\times G)/G$ be the projection. Then we define
a distribution $\bar H$ on $Ad\ P$ as follows:
$$\bar H_{[p,g]}=dq(H_p\oplus 0_q)$$
for any $[p,g]\in Ad\ P$. Invariance of $H$ with respect to the action of $G$ implies
that $\bar H$ is well-defined.

Also let $\bar V$ denote the `vertical' distribution defined by
$$\bar V_{[p,g]}=T_{[p,g]} (Ad\ P)_{\pi(p)}\subset T_{[p,g]} Ad\ P$$
for any $[p,g]\in Ad\ P$. Then we
have the decomposition of the tangent bundle of $Ad\ P$,
$$T Ad\ P =\bar V\oplus \bar H.$$
Therefore there is the projection $T Ad\ P \rightarrow \bar V$ coming from this
decomposition.

Any $p\in P$ is given, we let $\kappa_p$ now denote the homomorphism
$(Ad\ P)_{\pi(p)}\rightarrow G$ given by $\kappa_p[p, g]=g$
for any $g\in G$. Then if $h\in G$ we have $\kappa_{ph}=Ad_{h^{-1}}\kappa_p$.
Therefore, for any $\alpha\in \mathcal{H}^*(G)$, $\kappa_p^*\alpha\in \Omega^*((Ad\ P)_x; \mathbb{R})$
does not depend on the choice of $p\in P_x$, $x\in M$.
Now we may define, for any $\alpha \in \mathcal{H}^k(\alpha)$,
$\hat \alpha \in \Omega^k(Ad\ P;\mathbb{R})$ by the same formula as (6-1).
We again write
$$E_A: \mathcal{H}^*(G)\rightarrow \Omega^*(Ad\ P;\mathbb{R})$$
for the map
defined by $E_A(\alpha)=\hat\alpha$ for any $\alpha\in\mathcal{H}^*(G)$.

Again we observe:
\begin{lemma}
If $A$ is flat, $E_A(\alpha)\in \Omega^*(Ad\ P;\mathbb{R})$ is closed
for any $\alpha\in \mathcal{H}^*(G)$.
\end{lemma}
\bp
Assume $U\subset M$ is an open set for which a horizontal section
$s:U\rightarrow P_U$ exists. Then we consider the map
$\kappa_s: (Ad\ P)_U\rightarrow G$ defined by $\kappa_s[s(x),g]=g$
for any $x\in U$ and $g\in G$. The rest of the proof is a copy of the proof
of 6.1 above.
\ep

For any $p\in P_x$ we also have the inclusion $\iota_p:G\rightarrow (Ad\ P)_x\subset Ad\ P$
given by $\iota_p(g)=[p,g]$ for any $g\in G$, for which we have
$\iota_p^*\hat\alpha=\alpha$ for any $\alpha\in \mathcal{H}^*(G)$.

We have:
\begin{theorem}
Let $P\rightarrow M$, $G$ and $A$ be as in 6.2 above. Then there is the homomorphism of
algebras,
$$E_A:H^*_{\rm dR}(G)\rightarrow H^*_{\rm dR}(Ad\ P),$$
such that, for any $p\in P$,  $\iota_p^*E_A$
is the identity on $H^*_{\rm dR}(G)$.
\end{theorem}

Again we have:
\begin{corollary}
Let $P\rightarrow M$ and $G$ be as in 6.2 above. Then
there is an isomorphism of algebras:
$$H^*_{\rm dR}(Ad\ P)\cong H^*_{\rm dR}(M)\otimes H^*_{\rm dR}(G).$$
\end{corollary}

\section{The Maurer-Cartan 3-forms on the adjoint bundle} 
In this section we show that $E_A[\Theta]\in H^3_{\rm dR}(Ad\ P)$ does not depend on
the flat connection $A$, where
$E_A$ is the homomorphism introduced in the previous section
and $\Theta\in \Omega^3(G;\mathbb{R})$
is the Maurer-Cartan 3-form defined by a choice of an inner product
on $\mathcal{G}$. Furthermore we will show that the inner product can
be chosen so that $E_A[\Theta]$ is integral.

We begin by observing:
\begin{lemma}
Let $G$ be a compact connected semi-simple Lie group and
$P\rightarrow M$ be a flat $G$-bundle over a smooth manifold.
Then the homomorphism $E_A: H^3_{\rm dR}(G)\rightarrow H^3_{\rm dR}(Ad\ P)$ does not depend on the
choice of the flat connection $A$.
\end{lemma}
\bp
Recall that the Leray-Hirsch theorem applies to a flat bundle. Therefore
we have:
$$H^*_{\rm dR}(M)\otimes H^*_{\rm dR}(G)\cong H^*_{\rm dR}(Ad\ P)$$
where the isomorphism is given by $a\otimes b\rightarrow (\pi^*a)\wedge(E_A(b))$
for any $a\in H^*_{\rm dR}(M)$ and any $b\in H^*_{\rm dR}(G)$. Here $\pi$ denotes
the projection $Ad\ P\rightarrow M$.
In particular, since $H^i_{\rm dR}(G)=0$, for $i=1,2$,(see (3-6) above)
and $H^0_{\rm dR}(G)\cong \mathbb{R}$,
we have:
$$H^3_{\rm dR}(Ad\ P)=\pi^*H^3_{\rm dR}(M)\oplus E_A(H^3_{\rm dR}(G)).$$

Now consider the section $s_1: M\rightarrow Ad\ P$ which sends each $x\in M$ to the identity
element of $(Ad\ P)_x$. Note that $s_1^*\pi^*$ is the identity on $H^*_{\rm dR}(M)$ and
$s_1^* E_A:H^*_{\rm dR}(G)\rightarrow H^*_{\rm dR}(M)$ is the zero homomorphism. This shows that
$$Ker(H^3_{\rm dR}(Ad\ P)\overset{s_1^*}{\longrightarrow}H^3_{\rm dR}(M))=E_A(H^3_{\rm dR}(G)).$$
In particular, this means that $E_A(H^3_{\rm dR}(G))\subset H^3_{\rm dR}(Ad\ P)$ is independent
of the flat connection
$A$.

Now choose a $p\in P$ and observe that $\iota_p^*: E_A(H^3_{\rm dR}(G))\rightarrow H^3_{\rm dR}(G)$
is the inverse of $E_A: H^3_{\rm dR}(G)\rightarrow E_A( H^3_{\rm dR}(G))$ for any $A$. This proves
that $E_A: H^3_{\rm dR}(G)\rightarrow H^3_{\rm dR}(Ad\ P)$ does not depend on $A$.
\ep

It follows immediately that:
\begin{proposition}
Let $G$ and $P$ be as in 7.1 above and let $\Theta$ denote the Maurer-Cartan 3-form
of $G$ defined by choosing an Ad-invariant inner product on the Lie algebra of $G$.
Then the class $E_A[\Theta]\in H^3_{\rm dR}(Ad\ P)$ does not depend on the flat connection $A$.
\end{proposition}

When $G$ is compact connected and semi-simple,
now we may write $E: H^3_{\rm dR}(G)\rightarrow H^3_{\rm dR}(Ad\ P)$ to denote the map,
$E_A: H^3_{\rm dR}(G)\rightarrow H^3_{\rm dR}(Ad\ P)$, defined by choosing a flat connection
$A$. We have:

\begin{proposition}
Let $G$ and $P$ be as in 7.1 above.
Then there is a choice of Ad-invariant inner product on the Lie algebra of $G$
so that, if $\Theta$ is the Maurer-Cartan 3-form of $G$, $E[\Theta]\in H^3_{\rm dR}(Ad\ P)$ is
integral.
\end{proposition}
\bp
Consider the maps $\pi: Ad\ P\rightarrow M$ and $s_1:M\rightarrow Ad\ P$
introduced in the proof of 6.1 above
and the decomposition:
$$H^3_{\rm dR}(Ad\ P)=Im\ \pi^* \oplus Ker\ s_1^*.$$
We observed that $E( H^3_{\rm dR}(G))=Ker\ s_1^*$.

The decomposition above is valid also for the cohomology with
integer coefficient:
$$H^3(Ad\ P;\mathbb{Z})=Im\ \pi^* \oplus Ker\ s_1^*.$$
The homomorphism $j:H^3(Ad\ P;\mathbb{Z})\rightarrow H^3(Ad\ P;\mathbb{R})\equiv H^3_{\rm dR}(Ad\ P)$
also decomposes accordingly. Therefore we conclude the image of
$$Ker\ (H^3(Ad\ P;\mathbb{Z})\overset{s_1^*}{\longrightarrow}
H^3(M;\mathbb{Z}))\overset{j}{\longrightarrow} E( H^3_{\rm dR}(G))$$
is a full lattice in $E( H^3_{\rm dR}(G))$.

We have observed in \S 3 above that
the set $K$ of all $[\Theta]\in H^3_{\rm dR}(G)$ determined by choices of Ad-invariant
inner products on the Lie algebra of $G$ is an open cone in $H^3_{\rm dR}(G)$.
Therefore $E(K)$ is also an open cone in $E (H^3_{\rm dR}(G))$. We conclude
$$j(Ker\ s_1^*)\cap E(K)\neq \emptyset,$$
which proves the proposition.
\ep

Now let $G$ be any connected compact Lie group and $P\rightarrow M$
be any principal $G$-bundle. Let $Z_0$ denote the identity component of the center
of $G$. Then $G'=G/Z_0$ is a semi-simple Lie group and $P'=P/Z_0$ is a
principal $G'$-bundle. Write $q:G\rightarrow G'$ and $\hat q:P\rightarrow P'$
to denote the projections.

Let $A$ be a connection on $P$ and $H$ be the horizontal distribution of $P$
determined by $A$. Then $H$ determines a horizontal distribution $H'$ on $P'$:
Define $H'_{[p]}=d\hat q H_p$ for any $[p]=\hat q (p)\in P'$. Let $A'$ denote the
connection on $P'$ corresponding to $H'$.

Consider the map $\bar q :Ad\ P\rightarrow Ad\ P'$ defined by
$$\bar q[p,g]=[[p], [g]]$$
for any $(p,g)\in P\times G$. Let $\bar H$ and $\bar H'$ denote
the `horizontal' distributions respectively of $Ad\ P$ and $Ad\ P'$
determined respectively by $A$ and $A'$. Then we have that
$$d\bar q(\bar H_{[p,g]})=\bar H'_{[[p],[g]]}$$
for any $(p,g)\in P\times G$.

Let $\alpha\in\Omega^*(G';\mathbb{R})$ be a bi-invariant form.
Then $q^*\alpha\in \Omega^*(G;\mathbb{R})$
is also bi-invariant. Then we have the forms  $E_{A'}(\alpha)$ and $E_A(q^*\alpha)$
respectively on $Ad\ P'$ and on
$Ad\ P$.
Then we have that
\bdm
E_A(q^*\alpha)=\bar q^* E_{A'}(\alpha). \tag{7-1}
\edm
In fact this equation follows from the equality,
\bdm
dq d\kappa_p\pi_A=d\kappa_{[p]}\pi_{A'}d\bar{q}:T_{[p,g]}Ad\ P \rightarrow T_{[g]}G', \tag{7-2}
\edm
for any $(p,g)\in P\times G$.
The notations $\kappa_p, \kappa_{[p]}$ and $\pi_A, \pi_{A'}$ are as introduced in the previous
section. To see (7-2), given any $Y\in T_{[p,g]}Ad\ P$, one may choose a curve $[\gamma(t), \delta(t)]$
on $Ad\ P$ such that the velocity at $t=0$ is $Y$, where $\gamma$ is a horizontal curve in $P$
such that $\gamma(0)=p$ and $\delta$ is a curve in $G$ such that $\delta(0)=g$. Then it is straightforward
to see that both $dq d\kappa_p\pi_A(Y)$ and $d\kappa_{[p]}\pi_{A'}d\bar{q}(Y)$ are $(q\delta)^{\centerdot}(0)$.

Now we are ready to provide a proof of the following:
\begin{theorem}
Let $G$ be a compact connected Lie group with the Lie algebra $\mathcal{G}$ and
$P\rightarrow M$ be a flat
$G$-bundle over a manifold $M$. Let $\Theta$ denote the Maurer-Cartan
3-form of $G$ for some Ad-invariant inner product on $\mathcal{G}$.
Then $E_A[\Theta]\in H^3_{\rm dR}(Ad\ P)$ does not depend on the flat connection
$A$. Furthermore, there is a choice of an Ad-invariant inner product
on $\mathcal{G}$ so that $E_A[\Theta]$ is integral.
\end{theorem}
\bp
Let $G'$, $q:G\rightarrow G'$, $P'$, $\bar{q}:Ad\ P\rightarrow Ad\ P'$, etc., mean
the same as in the above.

Flatness of $A$ implies the same
for the induced connection $A'$ on $P'$:
If $s:U\rightarrow P$ is a horizontal lifting
of an open set $U\subset M$ with respect to $A$, $q's:U\rightarrow P'$ is also
such with respect $A'$.

Consider the decomposition $\mathcal{G}=\mathcal{Z}\oplus \mathcal{G}'$
where $\mathcal{Z}$ is the center
of $\mathcal{G}$ and $\mathcal{G}'=[\mathcal{G},\mathcal{G}]$.
Note that $\mathcal{G}'$ can be identified
with the Lie algebra of $G'$. Consider $\mathcal{G}'$ with the restriction of
the inner product of $\mathcal{G}$.
Then let $\Theta'$ denote the Maurer-Cartan 3-form of $G'$ determined by
this inner product on $\mathcal{G}'$. Then we have
$$q^*\Theta'=\Theta,$$
which is (3-8) above.

Apply (7-1) above to have
\bdm
E_A[\Theta]=\bar{q}^*E_{A'}[\Theta']. \tag{$\ast$}
\edm
Since $E_{A'}[\Theta']$ does not depend on $A'$ by Proposition 7.2 above,
we conclude that $E_A[\Theta]$
does not depend on $A$.

Now note that any Ad-invariant inner product on $\mathcal{G}'$ can be
extended to one on $\mathcal{G}$. Since there is an Ad-invariant inner product on $\mathcal{G}'$
so that $E_{A'}[\Theta']$ is integral by Proposition 7.3 above, equation $(\ast)$ above
proves also the second assertion of
the theorem.
\ep

\section{The Chern-Simons functional under gauge transformations} 
In this section we define $\deg \varphi$ exploiting the cohomology of
$Ad\ P$ so that (1-4) holds.
Our argument is based on the two technical lemmas below.

Let $P\rightarrow M$ be a flat bundle over a manifold with a structure group $G$.
Let $A\in \mathcal{A}(P)$ and let $cs(A)$ be the 3-form defined by
choosing a flat reference connection $A_0$.
Let $\varphi:P\rightarrow P$
be a gauge transformation and $u:P\rightarrow G$ be the map determined by $\varphi$.
Let $\theta$ and $\Theta$ respectively denote the Maurer-Cartan form and 3-form of $G$
and let $H_0$ denote the horizontal distribution given by $A_0$.
\begin{lemma}
We have:
$$cs(\varphi^*A)-cs(A)-(u^*\Theta)^{H_0}=d\langle (Ad_{u^{-1}}(A-A_0))\wedge(u^*\theta)^{H_0}\rangle.$$
\end{lemma}
\bp
We have that $\varphi^*A=A_0 +u^{-1}d_{A_0}u +Ad_{u^{-1}}(A-A_0)$. Also note that $u^{-1}du=u^*\theta$.
Therefore we may write: $\varphi^*A -A_0=(u^*\theta)^{H_0}+Ad_{u^{-1}}(A-A_0)$.

Then since $A_0$ is flat we have:
\bdm
\aligned
cs&(\varphi^*A)-cs(A)+\frac{1}{6}\langle(u^*\theta)^{H_0}\wedge[(u^*\theta)^{H_0}\wedge(u^*\theta)^{H_0})]\rangle\\
 &=\langle (u^*\theta)^{H_0}\wedge(Ad_{u^{-1}}F_A)\rangle -\frac{1}{2}\langle (u^*\theta)^{H_0}\wedge [Ad_{u^{-1}}(A-A_0)\wedge Ad_{u^{-1}}(A-A_0)]\rangle\\
 &\ \ \ \ \ \ \ \ \ \ \ \ \ \ \ \ \ -\frac{1}{2}\langle Ad_{u^{-1}}(A-A_0)\wedge [(u^*\theta)^{H_0}\wedge (u^*\theta)^{H_0}]\rangle
\endaligned
\edm

On the other hand the RHS of the equation in the lemma can be rewritten as
\bdm
\langle d_{A_0} (Ad_{u^{-1}}(A-A_0))\wedge(u^*\theta)^{H_0}\rangle
-\langle (Ad_{u^{-1}}(A-A_0))\wedge d_{A_0}((u^*\theta)^{H_0})\rangle. \tag{1}
\edm

Using the fact that $A_0$ is flat, we observe:

{\it Claim:} For any $\alpha\in \Omega^1(P; V)$ where $V$ is any real vector space
$$d_{A_0}(\alpha^{H_0})=(d\alpha)^{H_0}.$$

{\it Proof.} It suffices to consider horizontal vector fields $X, Y$. Observe
$$
\aligned
d_{A_0}(\alpha^{H_0})(X, Y)&= d(\alpha^{H_0})(X,Y)\\
                           &=X\alpha(Y)-Y\alpha(X)-\alpha^{H_0}([X,Y]).
\endaligned
$$
Since $H_0$ is integrable, $[X,Y]$ is horizontal. Thus we have
$$d_{A_0}(\alpha^{H_0})(X, Y)=d\alpha(X, Y),$$
which proves the claim.

\vspace{7pt}
Also recall that $d\theta=-\frac{1}{2}[\theta\wedge\theta]$. Then applying the claim above we have:
\bdm
d_{A_0}((u^*\theta)^{H_0})=(u^*d\theta)^{H_0}=-\frac{1}{2}[(u^*\theta)^{H_0}\wedge (u^*\theta)^{H_0}].
\tag{2}
\edm

Now we have(see the proof of 4.3 above):
\bdm
d_{A_0}(Ad_{u^{-1}}(A-A_0)) =(Ad_{u^{-1}}ad({u^{-1}}^*\theta))^{H_0}\wedge(A-A_0) +Ad_{u^{-1}}d_{A_0}(A-A_0).
\tag{3}
\edm

We rewrite the first term of the RHS of (3) as follows:
$$
\aligned
(Ad_{u^{-1}}ad({u^{-1}}^*\theta))^{H_0}\wedge(A-A_0)&=((Ad_{u^{-1}}ad({u^{-1}}^*\theta))\wedge(A-A_0))^{H_0}\\
&=Ad_{u^{-1}}[({u^{-1}}^*\theta)^{H_0}\wedge (A-A_0)].
\endaligned
$$
Since ${u^{-1}}^*\theta=-Ad_u u^*\theta$, we finally write:
\bdm
(Ad_{u^{-1}}ad({u^{-1}}^*\theta))^{H_0}\wedge(A-A_0)=-[(u^*\theta)^{H_0}\wedge Ad_{u^{-1}}(A-A_0)].
\tag{4}
\edm

Also we note that
\bdm
d_{A_0}(A-A_0)=d(A-A_0) +[A_0 \wedge (A-A_0)]=F_A -\frac{1}{2}[(A-A_0)\wedge(A-A_0)]. \tag{5}
\edm

We insert (4) and (5) to (3) to have:
\bdm
\aligned
d_{A_0}(Ad_{u^{-1}}(A-A_0))=-[(u^*\theta)^{H_0}\wedge & Ad_{u^{-1}}(A-A_0)] + Ad_{u^{-1}}F_A\\
 & -\frac{1}{2}[Ad_{u^{-1}}(A-A_0)\wedge Ad_{u^{-1}}(A-A_0)].
\endaligned \tag{6}
\edm
Then we insert (2) and (6) to (1) to conclude the lemma.
\ep

Let $P\rightarrow M$ and $G$ be as in the above.
Let $A$ be any connection on $P$ and $H$ be the horizontal distribution of $P$
given by $A$. Let $\varphi:P\rightarrow P$ and  $u: P\rightarrow G$
be as in the above and $\hat u :M\rightarrow Ad\ P$ be defined by $\hat{u}(x)=[p, u(p)]$
choosing any $p\in P_x$.
Let $\alpha\in \Omega^k(G; \mathbb{R})$ be a bi-invariant form. Then we have the form
$E_A(\alpha)\in \Omega^k(Ad\ P; \mathbb{R})$ as in \S 6 above. Also note that
$u^*\alpha$ is an invariant form on $P$ and therefore that $(u^*\alpha)^H$ can be regarded
as a form on $M$.

\begin{lemma}
We have the equality,
$(u^*\alpha)^H=\hat{u}^*(E_A(\alpha)),$
as forms on $M$.
\end{lemma}
\bp
Let $X_1,\cdots, X_k\in T_x M$. Let $\tilde X_1, \cdots, \tilde X_k\in T_pP$ be the
horizontal lifting with respect to $H$ where $p$ is an arbitrarily chosen point from $P_x$.
Then we have:
$$(u^*\alpha)^H(X_1,\cdots, X_k)=\alpha(du \tilde X_1, \cdots, du \tilde X_k).$$

On the other hand we have:
$$\hat{u}^*(E_A(\alpha))(X_1,\cdots, X_k)=(E_A(\alpha))(d\hat u X_1, \cdots, d\hat u X_k).$$
Recall from \S 6 the map $\kappa_p : (Ad\ P)_x \rightarrow G$ and
the projection $\pi_A: T(Ad\ P)\rightarrow \bar V$ given by the decomposition
$T(Ad\ P)= \bar V \oplus \bar H$.
Then we have:
$$(E_A\alpha)(d\hat u X_1, \cdots, d\hat u X_k)=\alpha(d\kappa_p \pi_A d\hat u X_1,  \cdots,d\kappa_p \pi_A d\hat u X_k).$$

The proof is complete if we observe:

{\it Claim:} For any $X\in T_x M$ let $\tilde X\in T_p P$ be the horizontal
lifting with respect to $A$. Then we have $d\kappa_p \pi_A d\hat{u} X= du \tilde X$.

{\it Proof:} Let $X$ be represented by a $C^1$-curve $\gamma:(-\epsilon,\epsilon)\rightarrow M$,
that is, assume $\gamma(0)=x$ and $\dot{\gamma}(0)=X$. Let
let $\tilde \gamma:(-\epsilon,\epsilon)\rightarrow P$ be the horizontal lifting such that $\tilde\gamma(0)=p$.
Then we have $\hat u(\gamma(t))=[\tilde\gamma(t), u(\tilde \gamma(t))]=q(\tilde\gamma(t), u(\tilde\gamma(t)))$,
where $q: P\times G\rightarrow Ad\ P$ is the projection. This means that we have:
$$d\hat u X=dq(\tilde X,du \tilde X)\in T_{[p,u(p)]}Ad\ P,$$
understanding $T_{(p,u(p))}(P\times G)\equiv T_p P\oplus T_{u(p)}G$.
It follows that  $\pi_A d\hat{u} X= dq(0, du \tilde X)$, which is
represented by the curve $[p, u(\tilde\gamma(t))]$. Therefore $d\kappa_p\pi_A d\hat{u} X$
is represented by the curve $u(\tilde\gamma(t))$, which proves the claim.
\ep

Now let $M$ be a closed oriented $3$-manifold and $P\rightarrow M$ admit a flat
connection $A_0$. If $CS:\mathcal{A}(P)\rightarrow \mathbb{R}$ is the Chern-Simons
functional defined by the reference connection $A_0$, Lemma 8.1 and 8.2 above imply
that
$$CS(\varphi^*A)-CS(A)=\int_M\hat{u}^*E_{A_0}(\Theta).$$
In particular, note that the 2-form $\langle (Ad_{u^{-1}}(A-A_0))\wedge(u^*\theta)^{H_0}\rangle$
in the equation of Lemma 8.1 is invariant and horizontal and therefore
the equation holds on $M$.

If $G$ is connected and compact, $E_{A_0}[\Theta]\in H^3_{\rm dR}(Ad\ P)$ does not depend on the flat
connection $A_0$ by Theorem 7.4 above. Therefore $\deg\varphi\in\mathbb{R}$
is well-defined by the following:
\begin{defn}
Let $P\rightarrow M$ be a flat bundle over an oriented closed $3$-manifold
whose structure group is connected and compact. If $\varphi:P\rightarrow P$ is a gauge transformation
and $\hat u:M\rightarrow Ad\ P$ is the map determined by $\varphi$, we define $\deg\varphi$ by
$$\deg\varphi = \int_M\hat u^*E_A(\Theta)\in \mathbb{R}$$
by choosing any flat connection $A$ on $P$.
\end{defn}

Furthermore, under the assumption that the structure group is compact and connected,
Theorem 7.4 above says that
there exists an Ad-invariant inner product on the Lie algebra so that
$E_A[\Theta]$ is an integral class of $H^3_{\rm dR}(Ad\ P)$ for any flat connection $A$.
Therefore, by adopting this
inner product, we may let $\deg\varphi$
be an integer for any gauge transformation $\varphi$, if necessary.

Since the choice of a reference connection affects the functional
$CS: \mathcal{A}(P)\rightarrow \mathbb{R}$ only by addition of a constant, by Lemma 8.1 and 2 above we have:
\begin{theorem}
Let $P\rightarrow M$ be a flat bundle over an oriented closed $3$-manifold
with a connected compact structure group. Let $CS$ denote the Chern-Simons functional
defined by choosing a reference connection. Then for any gauge transformation
$\varphi:P\rightarrow P$
and any connection $A$ on $P$ we have:
$$CS(\varphi^*A)-CS(A)=\deg \varphi .$$
\end{theorem}

\end{document}